\pdfoutput=1
\RequirePackage{ifpdf}
\ifpdf % We are running pdfTeX in pdf mode
\documentclass[pdftex]{sigma}
\else
\documentclass{sigma}
\fi

{ \theoremstyle{definition}
\newtheorem*{comments}{Comments}}

\begin{document}

\allowdisplaybreaks

\renewcommand{\PaperNumber}{078}

\FirstPageHeading

\ShortArticleName{Some Noncommutative Matrix Algebras Arising in the Bispectral Problem}

\ArticleName{Some Noncommutative Matrix Algebras\\
Arising in the Bispectral Problem}

\Author{F.~Alberto GR\"UNBAUM}

\AuthorNameForHeading{F.A.~Gr\"unbaum}

\Address{Department of Mathematics, University of California, Berkeley, CA 94720 USA}
\Email{\href{mailto:grunbaum@math.berkeley.edu}{grunbaum@math.berkeley.edu}}
\URLaddress{\url{http://math.berkeley.edu/~grunbaum/}}

\ArticleDates{Received May 01, 2014, in f\/inal form July 17, 2014; Published online July 24, 2014}

\Abstract{I revisit the so called ``bispectral problem'' introduced in a~joint paper with Hans Duistermaat a~long time
ago, allowing now for the dif\/ferential operators to have matrix coef\/f\/icients and for the eigenfunctions, and one of the
eigenvalues, to be matrix valued too.
In the last example we go beyond this and allow both eigenvalues to be matrix valued.}


\Classification{13N10; 16S32; 35P05}

\section{Introduction}

\looseness=-1
The characterization of certain {\it commutative} algebras of dif\/ferential operators goes back to pio\-nee\-ring work of
Schur and Burchnall jointly with Chaundy in the early part of the twenty century, see~\cite{BC1, BC2, BC3, Sch}.
The subject remained dormant until the 1970's where the study of the Korteweg--de Vries equation revived and enriched
the interaction between algebraic geometry, dif\/ferential operators, group representation theory and many other parts of
mathematics.
Suddenly the KdV equation, and related ones, such as the Toda and KP equations, started showing up in the solution of
many ``unrelated'' problems.
For an account of this revival see, for instance,~\cite{HK, K,M}.

One such problem was posed and solved in~\cite{DG} and has become known as the ``bispectral problem''.

It consists of f\/inding all scalar coef\/f\/icient ordinary dif\/ferential operators~$L$ and~$B$ in the variables~$x$ and~$z$
respectively such that there is a~common eigenfunction $\psi(x,z)$ satisfying the eigenvalue equations
\begin{gather*}
%\label{bispdef1}
L\psi =p(z)\psi,
\qquad
B\psi=\theta(x)\psi
\end{gather*}
for nonconstant scalar valued functions~$p$ and~$\theta$.

More explicitly we wanted to f\/ind all nontrivial instances where a~function $\psi(x,z)$ satisf\/ies
\begin{gather*}
L\left(x,\frac {d}{dx} \right)\psi(x,z) \equiv \left(- \left(\frac {d}{dx}\right)^2 + V(x)\right)\psi(x,z) = z\psi(x,z)
\end{gather*}
as well as
\begin{gather*}
B\left(z,\frac {d}{dz} \right)\psi(x,z) \equiv \left(\sum\limits_{i=0}^M b_i(z)\left(\frac {d}{dz} \right)^i\right)
\psi(x,z) = \theta(x)\psi(x,z).
\end{gather*}
All the scalar valued functions $V(x),b_i(z),\theta(x)$ are, in principle, arbitrary except for smoothness assumptions.
Notice that here~$M$ is arbitrary (f\/inite).
The operator~$L$ could be of higher order, but in~\cite{DG} we stick to order two.

The paper with Duistermaat starts by proving the importance of the so called ``ad-conditions'', see~\eqref{adcond}, an extremely
complicated system of non-linear ordinary dif\/ferential equations in the coef\/f\/icients of~$L$ and the function~$\theta$.
This line of attack, which produced the f\/irst few examples of bispectral situations, would have led nowhere if it were
not for the {\it totally unexpected} observation that the (rational solutions of the) KdV equation was lurking around.

Maybe the most remarkable fact is that the answer to this innocent looking question hides connections with many
important developments in the area of integrable systems.
Among these, one f\/inds a~useful role for the Darboux process, the notion of monodromy and f\/inally the appearance of
certain f\/lows, such as KP, that play a~central role.
The set of all possible~$L$ and the algebra of all dif\/ferential operators~$B$ going with a~given~$L$ was explicitly
characterized in~\cite{DG} in terms of an appropriate ``tau'' function.
In the case discussed in~\cite{DG} there are two families of potentials~$V(x)$ going with bundles of rank one and two
respectively.
In the f\/irst case the~``tau'' functions go back to Schur and one is dealing with the KdV f\/lows; in the second case,
called the even family in~\cite{DG}, one is dealing with the Virasoro f\/lows (the master symmetries of KdV).

The ``tau'' functions were given explicitly in~\cite{DG} and the observation regarding the second f\/low came later,
see~\cite{MZ}.
As far as I know no one has yet traced these second ``tau'' functions to something like characters of certain group
representation.

As we just mentioned, in the scalar valued case the {\it commutative algebra} of all dif\/ferential operators~$B$ going
with a~f\/ixed~$L$ is explicitly characterized in terms of an appropriate ``tau'' function.
This description of this {\it commutative algebra} is given below.
Notice that we also have a~characterization of the possible potentials~$V(x)$ that appear in the operator~$L$.
Nothing like this is available in the noncommutative case, all that we have at this point are a~few examples, as we
will see later.

For the $V(x)$ in the KdV family we have
\begin{gather*}
V(x) = \sum\limits_{p \in {\mathcal P}} \frac {\nu_p(\nu_p+1)}{(x-p)^2}
\end{gather*}
with ${\mathcal P}$ a~f\/inite subset of ${\mathbb C}$, and $\nu_p \in {\mathbb Z}_{>0}$ for $p \in {\mathcal P}$ being
such that
\begin{gather*}
\sum\limits_{\substack{q \in {\mathcal P}
\\
q \ne p}} \frac {\nu_q(\nu_q+1)}{(q-p)^{2j+1}} = 0 \qquad \mbox{for} \ \ 1 \le j \le \nu_p \quad \text{and each} \quad p \in {\mathcal P}.
\end{gather*}
One can also write
\begin{gather*}
V(x) = -2 \left(\frac {\theta'(x)}{\theta(x)} \right)'.
\end{gather*}
Here
\begin{gather*}
\theta(x) = \prod\limits_p (x-p)^{\frac {1}{2} \nu_p(\nu_p+1)},
\end{gather*}
and~$p$ runs over ${\mathcal P}$.

For these potentials we have the following characterization of the algebra of dif\/ferential operators in the spectral
variable:

{\it The eigenfunction $\psi(x,z)$ satisfies an equation of the form $B(z,\partial_z)\psi(x,z) = \theta(x)\psi(x,z)$ if
and only if the polynomial~$\theta$ has the property that}
\begin{gather*}
\theta^{(2j-1)}(p) = 0 \qquad \mbox{for all} \ \ 1 \le j \le \nu_p, \quad \mbox{for each pole} \ \ p \in {\mathbb C} \ \ \mbox{of} \ \ V.
\end{gather*}

For the $V(x)$ in the even family we have a~similar situation.

In~\cite{GH} ones sees how a~q-version (replacing dif\/ferential operators by $q$-dif\/ferences) of the scalar valued
bispectral problem touches on some very important work of Dick Askey.
This is one more instance of the wide impact of ``special functions'' in several parts of mathematics.

Several authors have made important contributions to illuminate this ``bispectral problem'' from dif\/ferent perspectives.
Since this paper is not intended to be a~survey I have restrained myself from mentioning most of this very nice work.
Apologies, apologies, \dots.
Instead of surveying all that is known in the {\it scalar} case I want to take a~jump back in time to the early 1980's
and pose some conjectures dealing with a~{\it matrix} valued version of the problem that Hans Duistermaat and I started
looking at back then.
In this much more dif\/f\/icult situation there is still no clear picture of a~complete solution.

As pointed out in~\cite{BGK} the f\/irst work on a~matrix valued version of the bispectral problem is due to J.P.~Zubelli (in his Berkeley Ph.D.~thesis), see for instance~\cite{Z1}.
For more recent work along these lines, see~\cite{GKS,SZ}.

The transition to a~matrix valued version of the same problem, with one of the two variables becoming discrete, is
featured for instance in a~series of papers with Pacharoni and Tirao, see~\cite{GPT1,GPT2}, involving matrix valued
spherical functions for certain symmetric spaces.
Here the spherical functions (properly ``packaged'') give rise to matrix valued orthogonal polynomials (a notion due to
M.G.~Krein~\cite{K2, K1}) that happen to satisfy dif\/ference as well as dif\/ferential equations.
Other examples of matrix valued orthogonal polynomials were given in joint work with Duran~\cite{DG1}.
In all these papers, starting with~\cite{GPT1} the observation is made that it is important that the operators~$\mathcal
L$ and~$\mathcal B$ should act on the left and the right respectively.
There is no other way to insure the commutativity of these two operators with matrix valued coef\/f\/icients.
This is also noticed in a~short paper with Iliev~\cite{GI}, where we try to use the ad-conditions.

The problem is now formulated in the form
\begin{gather}
\label{bispdef2}
\mathcal L\Psi =p(z)\Psi,
\qquad
\Psi\mathcal B=\Theta(x)\Psi
\end{gather}
for nonconstant functions~$p$ and~$\Theta$ which are scalar and {\it matrix valued} respectively.
The dif\/ferential operators have matrix valued coef\/f\/icients and they act on the matrix valued eigenfunction~$\Psi$.

The problem of classifying the noncommutative algebra of dif\/ferential operators $\mathcal B$ going with a~f\/ixed
$\mathcal L$, i.e., with a~f\/ixed family of matrix valued orthogonal polynomials, was f\/irst considered in a~joint paper
with Castro~\cite{CG}.
There we formulate conjectures about the structure of the algebra in a~few basic examples.
For a~matrix valued version of the Hermite polynomials the conjecture was eventually proved by Tirao~\cite{Tir} in
a~paper that gives a~complete ``nice presentation'' of this algebra in terms of generators and relations.

For a~f\/ixed $\mathcal L$ (going with a~f\/ixed family of matrix valued orthogonal polynomials) the algebra in question is
isomorphic to the algebra of all matrix valued polynomials $\Theta(x)$ satisfying the ``ad-conditions'' mentioned in the
introduction
\begin{gather}
\label{adcond}
(\operatorname{ad} \mathcal L)^{m+1}(\Theta(x))=0
\end{gather}
for some~$m$.
This is the starting point of my work with Duistermaat and has been seen to hold in the matrix case too,
see~\cite{GrTi}.

My original motivation for the bispectral problem came from a~very concrete application: the remarkable observation~by
D.~Slepian, H.~Landau and H.~Pollak at Bell Labs back in the 1960's that the integral ``time-and-band limiting''
operator of Claude Shannon allowed for an explicit commuting dif\/ferential operator, see~\cite{G2,G1,S2,S1}.

In a~few cases I have been able to go back from bispectral situations to the motivating problem,
see~\cite{G3,G5,G4,GLP}, see also~\cite{P1,P2}.
Very recently this has been accomplished for the f\/irst time in a~noncommutative set-up, see~\cite{GPZ}.

The work in~\cite{GLP} has been picked up in the applied literature, see for instance~\cite{SD,SDW}.

\section{Contents of the paper}

To the best of my knowledge the f\/irst papers to consider the case of two dif\/ferential operators acting from dif\/ferent
directions are~\cite{BL} and~\cite{BGK}.
These papers acknowledge that this idea had been advanced in the papers mentioned above dealing with dif\/ferential and
dif\/ference operators.

The f\/irst paper~\cite{BL} includes two examples: in each case $\mathcal L$ has order two and $\mathcal B$ has order
three in one case and four in the other one.
Here the authors show the bispectrality of two examples given in~\cite{GV} of (matrix) Schr\"odinger operators with
trivial monodromy.
Note that in the abstract of~\cite{BL} both eigenvalues are allowed to be {\em matrix valued}.
In the actual examples in~\cite{BL} they are both scalar valued and in the body of the paper (see the proof of Theorem~3.6)
the right hand side of the equation involving~$\mathcal B$ is written dif\/ferently that the way it appears in the
abstract.
In a~recent conversation with J.~Liberati he seems to prefer the second version.
I~prefer the other one, the one in their abstract.

The second paper~\cite{BGK} exhibits an example arising from spin Calogero systems and the use of the Wilson
bispectral involution, a~remarkably useful tool, see~\cite{W4, W1,W2,W3}.
In~\cite{BGK} the opera\-tor~$\mathcal L$ has order four and $\mathcal B$ has order six.
The problem is formulated there with both eigenvalues taken as scalar valued functions.
Notice that in some of the examples of matrix bispectrality that arise in~\cite{CG,DG1,GI,GPT1,GPT2} as well as others
in the context of matrix valued orthogonal polynomials this is not necessarily the case, and allowing (at least) one of
the eigenvalues to be matrix valued is the natural thing to do.

I f\/irst deal with the f\/irst example in~\cite{BL}, see also~\cite{GV}, where a~bispectral situation is given involving
a~matrix valued eigenfunction but restricted to the case of scalar valued eigenvalues.
I do allow for arbitrary matrix valued polynomial $\Theta(x)$.

I then consider an analog of this two-by-two example in a~three-by-three step-up.
In each case I give an explicit conjecture about the algebra of all possible $\Theta(x)$, or equivalently all possible
operators $\mathcal B$.

Finally I give an algebra and a~conjecture about all possible matrix valued eigenvalues related to an example which is
connected to spin Calogero models but more elaborate than the one in~\cite{BGK}.

\section{The f\/irst example}

Take for $\Psi(x,z)$ the matrix valued function
\begin{gather*}
\Psi(x,z) = e^{xz}
\begin{pmatrix}
z - x^{-1} & x^{-2}
\\
0 & z - x^{-1}
\end{pmatrix}
\end{gather*}
and consider all instances of matrix valued polynomials $\Theta(x)$ and dif\/ferential operators $\mathcal B$ (with matrix
coef\/f\/icients $b_i(z)$) such that
\begin{gather*}
\Psi \mathcal B \equiv \sum\limits_{i=1}^m \big(\partial_z^i\Psi\big) b_i = \Theta(x)\Psi(x).
\end{gather*}
In this case one has
\begin{gather*}
\mathcal L\Psi = -z^2\Psi
\end{gather*}
with
\begin{gather*}
\mathcal L = -\partial_x^2 + 2
\begin{pmatrix}
x^{-2} & -2 x^{-3}
\\
0 & x^{-2}
\end{pmatrix}.
\end{gather*}
In other words for this specif\/ic dif\/ferential operator in the variable~$x$ we are asking for all bispectral ``partners''
of~$\mathcal L$.

In~\cite{BL} one f\/inds that one such pair $(\mathcal B,\Theta)$ is given~by
\begin{gather*}
\mathcal B = \partial_z^3 - 3\partial_z^2 \frac {1}{z} + 3\partial_z \frac {1}{z^2} + 3
\begin{pmatrix}
0 & z^{-2}
\\
0 & 0
\end{pmatrix}
\end{gather*}
and $\Theta(x)$ the scalar-valued polynomial
\begin{gather*}
\Theta(x) = x^3.
\end{gather*}

For operators $\mathcal L$ of the form
\begin{gather*}
\mathcal L = -\partial_x^2 + \mathbb U(x)
\end{gather*}
with a~matrix valued potential $\mathbb U(x)$ one can argue as in~\cite{DG} (using the ad-conditions) and conclude that
$\Theta(x)$ has to be a~polynomial with matrix valued coef\/f\/icients.
The ``matrix valued'' eigenvalue $\Theta(x)$ has to satisfy lots of other restrictions for the equation
\begin{gather*}
\Psi \mathcal B = \Theta(x)\Psi
\end{gather*}
to hold, and the problem at hand is to describe explicitly all these $\Theta(x)$.
An important observation is that this set of~$\Theta$'s form a~noncommutative algebra of polynomials in~$x$ and the
algebra of the corresponding dif\/ferential operators $\mathcal B$ is isomorphic to it.

\begin{conjecture} In our case the set of all~$\Theta$'s is the algebra of all polynomials of the form
\begin{gather*}
\begin{pmatrix}
r^{11}_0 & r^{12}_0
\vspace{1pt}\\
0 & r^{11}_0
\end{pmatrix}
+
\begin{pmatrix}
r^{11}_1 & r^{12}_1
\vspace{1pt}\\
0 & r^{11}_1
\end{pmatrix}
  x +
\begin{pmatrix}
r^{11}_2 & r^{12}_2
\vspace{1pt}\\
r^{11}_1 & r^{22}_2
\end{pmatrix}
  x^2 +
\begin{pmatrix}
r^{11}_3 & r^{12}_3
\vspace{1pt}\\
r^{22}_2 + r^{11}_2 - r^{12}_1 & r^{22}_3
\end{pmatrix}
  x^3 + x^4P(x)  ,
\end{gather*}
where $P(x)$ is an arbitrary $2 \times 2$ matrix valued polynomial and all the variables
$r^{11}_0$, $r^{12}_0$, $r^{11}_1$, $r^{12}_1$, $r^{11}_2$, $r^{12}_2$, $r^{22}_2$, $r^{11}_3$, $r^{12}_3$, $r^{22}_3$ are arbitrary.

Moreover, for each such~$\Theta$ one can give an explicit expression for the corresponding operator~$\mathcal B$.
\end{conjecture}

\begin{comments} It is not hard to check that the set of polynomials $\Theta(x)$ given above forms an algebra.
It is also not hard to propose a~list of generators for the algebra and certain relations among them.
The remaining open problem is to prove the conjecture and to give a~``nice description'' of the algebra in terms of
generators and relations.

The situation is analogous to what was put forth in~\cite{CG}, namely a~conjecture about the algebra and a~collection of
generators for a~few examples.
The proof of these conjectures and a~``nice description'' of the algebra was only done, for one of the examples put
forward in~\cite{CG} and~\cite{Tir}.
\end{comments}

\section{Second example}

Take for $\Psi(x,z)$ the matrix valued function
\begin{gather*}
\Psi(x,z) = \left[\partial_x -
\begin{pmatrix}
x^{-1} & -x^{-2} & x^{-3}
\\
0 & x^{-1} & -x^{-2}
\\
0 & 0 & x^{-1}
\end{pmatrix}
\right] e^{xz}I = e^{xz}
\begin{pmatrix}
z - x^{-1} & x^{-2} & - x^{-3}
\\
0 & z - x^{-1} & x^{-2}
\\
0 & 0 & z - x^{-1}
\end{pmatrix}
.
\end{gather*}
Here one can see that
\begin{gather*}
\mathcal L\Psi = -z^2\Psi
\end{gather*}
with
\begin{gather*}
\mathcal L = -\partial_x^2 + 2
\begin{pmatrix}
x^{-2} & -2 x^{-3} & 3 x^{-4}
\\
0 & x^{-2} & -2x^{-3}
\\
0 & 0 & x^{-2}
\end{pmatrix}.
\end{gather*}

\begin{conjecture} The algebra of all matrix valued polynomials $\Theta(x)$ for which there exist some operator
$\mathcal B$ with
\begin{gather*}
\Psi \mathcal B = \Theta(x)\Psi
\end{gather*}
is the algebra of all polynomials of the form
\begin{gather*}
\begin{pmatrix}
r^{11}_0 & r^{12}_0 & r^{13}_0
\vspace{1pt}\\
0 & r^{22}_0 & r^{23}_0
\vspace{1pt}\\
0 & 0 & r^{11}_0
\end{pmatrix}
+
\begin{pmatrix}
r^{11}_1 & r^{12}_1 & r^{13}_1
\vspace{1pt}\\
r^{22}_0 - r^{11}_0 & r^{22}_1 & r^{23}_1
\vspace{1pt}\\
0 & r^{22}_0 - r^{11}_0 & r^{11}_1 + r^{23}_0 - r^{12}_0
\end{pmatrix}
 x
\\
\qquad{}
+
\begin{pmatrix}
r^{11}_2 & r^{12}_2 & r^{13}_2
\vspace{1pt}\\
r^{22}_1 - r^{11}_1 - r^{23}_0 + r^{12}_0 & r^{22}_2 & r^{23}_2
\vspace{1pt}\\
r^{22}_0 - r^{11}_0 & r^{22}_1 - r^{11}_1 & r^{11}_2 + r^{23}_1 - r^{12}_1
\end{pmatrix}
  x^2
\\
\qquad{}+
\begin{pmatrix}
r^{11}_3 & r^{12}_3 & r^{13}_3
\vspace{1pt}\\
r^{21}_3 & r^{22}_3 & r^{23}_3
\vspace{1pt}\\
r^{22}_1 - 2r^{11}_1 - r^{23}_0 + r^{12}_0 & r^{32}_3 & r^{33}_3
\end{pmatrix}
  x^3  \\
\qquad{}
+
\begin{pmatrix}
r^{11}_4 & r^{12}_4 & r^{13}_4
\vspace{1pt}\\
r^{21}_4 & r^{22}_4 & r^{23}_4
\vspace{1pt}\\
r^{32}_3 + r^{21}_3 - r^{22}_2 - r^{11}_2 + r^{12}_1 & r^{32}_4 & r^{33}_4
\end{pmatrix}
  x^4
\\
\qquad{}
+
\begin{pmatrix}
r^{11}_5 & r^{12}_5 & r^{13}_5
\vspace{1pt}\\
r^{21}_5 & r^{22}_5 & r^{23}_5
\vspace{1pt}\\
r^{32}_4 + r^{21}_4 - r^{33}_3 - r^{22}_3 - r^{11}_3 + r^{23}_2 + r^{12}_2 - r^{13}_1 & r^{32}_5 & r^{33}_5
\end{pmatrix}
  x^5 + x^6 P(x),
\end{gather*}
where $P(x)$ is an arbitrary $3 \times 3$ matrix valued polynomial and all the variables
$r^{11}_0,r^{12}_0,\dots,r^{33}_5$ are arbitrary.
\end{conjecture}

Once again one can prove by an explicit computation that this set of~$\Theta$'s forms an algebra and one can single out
generators and some relations among them.
What is needed is a~proof of the conjecture and a~``nice presentation'' in terms of generation and relations.
Here is a~very simple case of the general result above.
We can check that
\begin{gather*}
\frac {\partial^2\Psi}{\partial z^2}
\begin{pmatrix}
0 & 0 & 0
\\
0 & 0 & 0
\\
1 & 0 & 0
\end{pmatrix}
 + \frac {\partial \Psi}{\partial z}
\begin{pmatrix}
0 & 0 & 0
\\
1 & 0 & 0
\\
-2z^{-1} & 1 & 0
\end{pmatrix}
 + \Psi
\begin{pmatrix}
1 & 0 & 0
\\
-2{z}^{-1} & 2 & 0
\\
0 & 0 & 1
\end{pmatrix}
=
\begin{pmatrix}
1 & 0 & 0
\\
x & 2 & 0
\\
x^2 & x & 1
\end{pmatrix}
\Psi.
\end{gather*}

\section{Third example}

We consider now a~situation that will be explained fully in~\cite{GJZ}.
It comes about by looking at examples linked to the spin Calogero systems discussed for instance in~\cite{BGK}.

Consider the function $\Psi(x,z)$ given~by
\begin{gather*}
\Psi(x,z) = \frac {e^{xz}}{(x-2)xz} \left(
\begin{matrix} \dfrac {x^3z^2-2x^2z^2-2x^2z+3xz+2x-2}{xz} & \dfrac {1}{x} \vspace{2mm}
\\
\dfrac {(xz-2)}{z} & x^2z-2xz-x+1
\end{matrix}
\right),
\end{gather*}
which satisf\/ies
\begin{gather*}
\mathcal L\Psi = \Psi F(z)
\end{gather*}
with
\begin{gather*}
\mathcal L \equiv \left(
\begin{matrix} 0 & 0   \\
0 & 1
\end{matrix}
\right) \partial_x^2 + \left(
\begin{matrix} 0 & \dfrac {1}{(x-2)x^2} \vspace{2mm}
\\
-\dfrac {1}{x-2} & 0
\end{matrix}
\right) \partial_x + \left(
\begin{matrix} -\dfrac {1}{x^2(x-2)^2} & \dfrac {x-1}{x^3(x-2)^2} \vspace{2mm}
\\
\dfrac {2x-1}{x(x-2)^2} & -\dfrac {2x^2-4x+3}{x^2(x-2)^2}
\end{matrix}
\right)
\end{gather*}
and
\begin{gather*}
F(z) = \left(
\begin{matrix} 0 & 0 \\
0 & z^2
\end{matrix}
\right).
\end{gather*}
We can now check that
\begin{gather*}
\frac {\partial^2\Psi}{\partial z^2} \left(
\begin{matrix} 0 & 0 \vspace{2pt}   \\
-\dfrac {2z+1}{z} & 0
\end{matrix}
\right)  + \frac {\partial \Psi}{\partial z} \left(
\begin{matrix} 1 & 0 \vspace{2pt}
\\
\dfrac {2(z-1)}{z^2} & 1
\end{matrix}
\right)
 + \Psi \left(
\begin{matrix} -z^{-1} & 0 \vspace{2pt}
\\
 6 z^{-3} &  z^{-1}
\end{matrix}
\right) = \left(
\begin{matrix} x & 0
\\
x^2(x-2) & x
\end{matrix}
\right) \Psi.
\end{gather*}

We are thus dealing with a~situation of the following kind
\begin{gather*}
%\label{bispdef3}
\mathcal L\Psi =\Psi F(z),
\qquad
\Psi\mathcal B=\Theta(x)\Psi
\end{gather*}
for nonconstant matrix valued functions~$F$ and~$\Theta$.
As in~\eqref{bispdef2} the dif\/ferential operators have matrix valued coef\/f\/icients and they act on the matrix valued
eigenfunction~$\Psi$.

One can see that, for the $\Psi(x,z)$ above, the algebra of all $F(z)$ such that for some $\mathcal L$ one has
\begin{gather*}
\mathcal L\Psi = \Psi F(z)
\end{gather*}
is given by polynomials in~$z$ of the form
\begin{gather*}
\begin{pmatrix}
a~& 0
\\
b-a & b
\end{pmatrix}
+
\begin{pmatrix}
c & c
\\
a~-- b -c & -c
\end{pmatrix}
  z +
\begin{pmatrix}
a-b-c & c+a-b
\\
d & e
\end{pmatrix}
  \big(z^2\big)/2 + z^3 P(z),
\end{gather*}
where $P(z)$ is an arbitrary $2 \times 2$ matrix valued polynomial and all the variables $a$, $b$, $c$, $d$, $e$ are arbitrary.

It is not hard to see that this forms an algebra for which, once again, a~nice description in terms of generators and
relations remains a~challenge.

\subsection*{Acknowledgements}

The author is extremely grateful to referees who did an outstanding job in suggesting improvements to an earlier version
of this paper.

\pdfbookmark[1]{References}{ref}

\LastPageEnding

\end{document}